\documentclass[11pt]{article}

\usepackage{amsmath,amsfonts,amssymb,longtable}

\usepackage{makeidx}
\usepackage{amssymb}
\usepackage{amsfonts}
\usepackage{amsmath}
\usepackage{euscript}
\usepackage{theorem}
\usepackage{enumerate}

\theorembodyfont{\rm}      

\theoremstyle{change}      

\begingroup

\newtheorem{thm}{Theorem\hskip 5mm}[section]

\newtheorem{theorem}{Theorem\hskip 5mm}[section]

\newtheorem{note}[thm]{Note\hskip 5mm}

\endgroup

\begingroup

\endgroup

\begingroup

\endgroup

\newcommand{\qed}{{\unskip\nobreak\hfil\penalty 50\hskip 2em\hbox{}

 \nobreak\hfil$\square$\parfillskip=0pt\finalhyphendemerits=0\par}}

\def\GL{{\rm GL}}
\def\ki{{\chi(F)}}

\begin{document}

\begin{center}
{\bf EQUIVALENCE AND CONGRUENCE OF MATRICES UNDER THE ACTION OF
STANDARD PARABOLIC SUBGROUPS}
\end{center}

\begin{center}
{\sc FERNANDO SZECHTMAN}
\end{center}

\begin{center}
{\small\textit{Department of Mathematics \& Statistics, University
of Regina, Saskatchewan, Canada, S4S 0A2}}
\end{center}

\section{Introduction} \label{ch0}

We fix throughout a field $F$ and a positive integer $n\geq 2$.
Let $M$ stand for the space of all $n\times n$ matrices over $F$
and write $G=\GL_n(F)$ for the general linear group.

We denote by $X'$ the transpose of $X\in M$. If $H$ is a subgroup
of $G$ and $X,Y\in M$ we say that $X$ and $Y$ are $H$-{\em
equivalent} if there exist $h,k\in H$ such that $Y=h'Xk$, and
$H$-{\em congruent} if there exists $h\in H$ such that $h'Xh=Y$.

The goal of this paper is to find necessary and sufficient
conditions for $H$-equivalence of arbitrary matrices, and
$H$-congruence of symmetric and alternating matrices, for various
subgroups $H$ of $G$, specifically the subgroups $U$, $B$ and
$P$, as defined below.

By $B$ we mean the group of all invertible upper triangular
matrices and by $U$ the group of all upper triangular matrices
whose diagonal entries are equal to 1.

We write $P$ for a standard parabolic subgroup of $G$, i.e. a
subgroup of $G$ containing $B$. Sections 8.2 and 8.3 of [3] ensure
that $P$ is generated by $B$ and a set $J$ of transpositions
(viewed as permutation matrices) of the form $(i,i+1)$, where
$1\leq i<n$. Let $e_1,...,e_n$ stand for the canonical basis of
the column space $F^n$. Consider the sequence of subspaces
$$
(0)\subset \langle e_1\rangle\subset \langle
e_1,e_2\rangle\subset\cdots\subset \langle
e_1,e_2,...,e_{n-1}\rangle\subset F^n,
$$
and let ${\cal C}$ be the chain obtained by deleting the $i$-th
intermediate term from the above chain if and only if $(i,i+1)$ is
in $J$. An alternative description for $P$ is that it consists of
all matrices in $G$ stabilizing each subspace in the chain ${\cal
C}$. Thus $P$ consists of block upper triangular matrices, where
each diagonal block is square and invertible. If ${\cal C}$ has
length $m$ then $2\leq m\leq n+1$ and the matrices in $P$ have
$m-1$ diagonal blocks, where the size of block $i$ is the
codimension of the $(i-1)$-th term of ${\cal C}$ in the $i$-th
term of ${\cal C}$, $1<i\leq m$.

In particular, if $m=2$ then $P=G$, while if $m=n+1$ then $P=B$.

We know that $G$-equivalence has the same meaning as rank
equality. It is also known that alternating matrices are
$G$-congruent if and only if they have the same rank. The same is
true for symmetric matrices under the assumptions that $F=F^2$
(every element of $F$ is a square) and $\ki\neq 2$ (the
characteristic of $F$ is not 2). If $F=F^2$ but $\ki=2$ then a
symmetric matrix is either alternating or $G$-congruent to a
diagonal matrix, and in both cases rank equality means the same as
$G$-congruence (see [8], chapter 1).

By a $(0,1)$-matrix we mean a matrix whose entries are all equal
to 0 or 1. A {\em sub-permutation} is a matrix having at most one
non-zero entry in every row and column. In exercise 3 of section
3.5 of [7] we find that every matrix $X$ is $B$-equivalent to a
sub-permutation $(0,1)$-matrix $Y$. If $X$ is invertible then the
uniqueness of the Bruhat decomposition yields that $Y$ is unique.
If $X$ is not invertible $Y$ is still unique, although it is
difficult to find a specific reference to this known result. Thus
two matrices are $B$-equivalent if and only if they share the same
associated sub-permutation $(0,1)$-matrix.

We may derive the problems of $B$-congruence and $U$-congruence
from their equivalence counterparts, except for symmetric matrices
in characteristic 2 when the situation becomes decidedly harder,
even under the assumption that $F=F^2$. One of our contributions
is a list of orbit representatives of symmetric matrices under $B$
and $U$-congruence if $\ki=2$ and $F=F^2$.

A second contribution is in regards to the question of
$P$-equivalence and $P$-congruence. We first determine various
conditions logically equivalent to $P$-equivalence. Let $W$ stand
for the Weyl group of $P$, i.e. the subgroup of $S_n$ generated by
$J$. We view $W$ as a subgroup of $G$. One of our criteria states
that two matrices $Y$ and $Z$ are $P$-equivalent if and only if
their associated sub-permutation $(0,1)$-matrices $C$ and $D$ are
$W$-equivalent. This generalizes the above results for $G$ and
$B$-equivalence. We also determine two alternative
characterizations of $P$-equivalence in terms of numerical
invariants of the top left block submatrices of $Y$ and $Z$, and
also of $C$ and $D$ (the well-known criterion for LU-factorization
using principal minors becomes a particular case).

Finally we show that for symmetric matrices (when $F=F^2$ and
$\ki\neq 2$) and alternating matrices, $P$-congruence has exactly
the same meaning as $P$-equivalence. We also furnish an
alternative characterization of $P$-congruence in terms of
$W$-conjugacy.

A restricted case of $P$-equivalence was considered is [5], but
our main results and goals are very distant from theirs. A
combinatorial study of $B$-congruence of symmetric complex
matrices is made in [1].

We remark that the congruence actions of $U$ on symmetric and
alternating matrices appear naturally in the study of a $p$-Sylow
subgroup $Q$ of the symplectic group $\mathrm{Sp}_{2n}(q)$ and the
special orthogonal group $\mathrm{SO}^+_{2n}(q)$, respectively.
Here $q$ stands for a power of the prime $p$. An investigation of
these actions is required in order to analyze the complex
irreducible characters of $Q$ via Clifford theory. We refer the
reader to [6] for details.

Suppose that $\ki\neq 2$ and $F=F^2$. At the very end of the paper
we count the number of orbits of symmetric matrices under
$B$-congruence. This number being finite, so is the number of
orbits of invertible matrices under $P$-congruence. We may
interpret this as saying that the double coset space $O\diagdown
G\diagup P$ is finite, where $O$ stands the orthogonal group. The
finiteness or not of a double coset space of the form $H\diagdown
G\diagup P$ for groups $G$ more general than ours has been studied
extensively.  Precise information can be found in the works of
Brundan, Duckworth, Springer and Lawther cited in the
bibliography.

We keep the above notation and adopt the following conventions.  A
$(1,-1)$-matrix is a sub-permutation alternating matrix whose only
non-zero entries above the main diagonal are equal to 1.

If $X\in M$ is a sub-permutation then there exists an $X$-{\em
couple} associated to it, namely a pair $(f,\sigma)$ where
$\sigma\in S_n$ and $f:\{1,...,n\}\to F$ is a function such that
$$
Xe_i=f(i)e_{\sigma(i)},\quad 1\leq i\leq n.
$$
We write $S(f)$ for the support of $f$, i.e. the set of points
where $f$ does not vanish.

\section{Equivalence Representatives under $U$ and $B$}

The Bruhat decomposition of $G$ can be interpreted as saying that
permutation matrices are representatives for the orbits of $G$
under $B$-equivalence. This can be pushed further by noting that
every matrix in $M$ is $B$-equivalent to a unique sub-permutation
$(0,1)$-matrix. We include a proof of this known result, which is
a particular case of Theorem \ref{main} below.

\begin{theorem}
\label{zero} Let $X\in M$. Then

(a) $X$ is $B$-equivalent to a unique sub-permutation
$(0,1)$-matrix.

(b) $X$ is $U$-equivalent to a unique sub-permutation matrix.
\end{theorem}

\noindent{\it Proof.} Existence is a simple exercise that we omit.

To prove uniqueness in (a) suppose that $Y$ and $Z$ are
sub-permutation $(0,1)$-matrices and that $c'Yd=Z$ for some
$c,d\in B$. Set $a=c'$ and $b=d^{-1}$, so that $aY=Zb$. We wish to
show that $Y=Z$.

Let $(f,\sigma)$ be a $Z$-couple and let $(g,\tau)$ be a
$Y$-couple, where $f,g:\{1,...,n\}\to \{0,1\}$. Notice that $S(f)$
and $S(g)$ have the same cardinality: the common rank of $Y$ and
$Z$.

We need to show that $S(f)=S(g)$ and that $\sigma(i)=\tau(i)$ for
every $i\in S(f)$.

As $a$ is lower triangular and $b$ is upper triangular, for all
$1\leq i\leq n$ we have
\begin{equation}
\label{e1}
Zbe_i=Z[b_{1i}e_1+\cdots+b_{ii}e_i]=b_{1i}f(1)e_{\sigma(1)}+\cdots+b_{ii}f(i)e_{\sigma(i)}
\end{equation}
and
\begin{equation}
\label{e2}
aYe_i=a[g(i)e_{\tau(i)}]=g(i)[a_{\tau(i),\tau(i)}e_{\tau(i)}+\cdots+a_{n,\tau(i)}e_{n}].
\end{equation}
Note that every diagonal entry of $a$ and $b$ must be non-zero.

Suppose that $i\in S(f)$. Then $e_{\sigma(i)}$ appears with
non-zero coefficient in (\ref{e1}), so it must likewise appear in
(\ref{e2}). We deduce that $i\in S(g)$ and $\tau(i)\leq
\sigma(i)$. This proves that $S(f)$ is included in $S(g)$. As they
have the same cardinality, they must be equal. Thus, for every $i$
in the common support of $f$ and $g$, we have $\tau(i)\leq
\sigma(i)$.

Suppose $\tau$ and $\sigma$ do not agree on $S(f)=S(g)$. Let $i$
be the first index in the common support such that
$\tau(i)<\sigma(i)$. Now $e_{\tau(i)}$ appears in (\ref{e2}) with
non-zero coefficient, so it must likewise appear in (\ref{e1}).
Thus we must have $\tau(i)=\sigma(j)$ for some $j$ such that $j<i$
and $j\in S(f)$. But then $\tau(j)=\sigma(j)=\tau(i)$, which
cannot be. This proves uniqueness in (a).

We use the same proof in (b), only that every diagonal entry of
$a$ and $b$ is now equal to 1, while $f$ and $g$ take values in
$F$. Then the old proof gives the additional information that
$f(i)=g(i)$ for all $i$ in the common support of $f$ and $g$, as
required. \qed

\section{Congruence Representatives under $U$ and $B$. Case 1}

\begin{thm}
\label{uno} Let $X\in M$.

(1) Suppose $\ki\neq 2$. If $X$ is symmetric then $X$ is
$U$-congruent to a unique sub-permutation matrix. Two symmetric
matrices are $U$-congruent if and only if they are $U$-equivalent.

(2) If $X$ is alternating then $X$ is $U$-congruent to a unique
sub-permutation matrix. Two alternating matrices are $U$-congruent
if and only if they are $U$-equivalent.
\end{thm}

\noindent{\it Proof.} Existence in (1) and (2) is a simple
exercise that we omit. Uniqueness in (1) and (2) follows from
uniqueness in Theorem \ref{zero}. Suppose $C$ and $D$ are
symmetric (resp. alternating) and $U$-equivalent. Let $Y,Z$ be
sub-permutation matrices $U$-congruent to $C$ and $D$,
respectively. Then $Y,Z$ are $U$-equivalent, so $Y=Z$ by Theorem
\ref{zero}. Hence $C$ and $D$ are $U$-congruent. The converse is
obvious.\qed

Much as above, we obtain the following result.

\begin{thm}
\label{Buno} Let $X\in M$.

(1) Assume $\ki\neq 2$ and $F=F^2$. If $X$ is symmetric then $X$
is $B$-congruent to a unique sub-permutation $(1,0)$-matrix. Two
symmetric matrices are $B$-congruent if and only if they are
$B$-equivalent.

(2) If $X$ is alternating then $X$ is $B$-congruent to a unique
$(1,-1)$-matrix. Two alternating matrices are $B$-congruent if and
only if they are $B$-equivalent.
\end{thm}

\section{Congruence Representatives under $U$ and $B$. Case 2}

We declare $X\in M$ to be a {\em pseudo-permutation} if $X$ is
symmetric, every column of $X$ has at most two non-zero entries,
and if there exists $j$ such that column $j$ of $X$ has two
non-zero entries then these must be $X_{jj}$ and $X_{ij}$ for some
$i<j$.

Suppose that $X$ is a pseudo-permutation matrix. Every pair
$(i,j)$ where $X_{ij}$ and $X_{jj}$ are non-zero is called an
$X$-{\em pair} (notice that $X_{ii}=0$ in this case). Suppose that
$(i,j)$ is an $X$-pair. If $(k,\ell)$ is also an $X$-pair we say
that $(k,\ell)$ is {\em inside} $(i,j)$ provided $i<k<\ell<j$. By
an $X$-index we mean an index $s$ such that $X_{ss}\neq 0$ and
this the only non-zero entry in column $s$ of $X$. If $s$ is an
$X$-index then $s$ is $X$-{\em interior} to the $X$-pair $(i,j)$
if $i<s<j$. An $X$-pair is {\em problematic} if it has an $X$-pair
inside it or an $X$-index interior to it.

We refer to $X$ as a {\em specialized pseudo-permutation} if $X$
is a pseudo-permutation with no problematic $X$-pairs.

As an illustration, the $(0,1)$-matrices
$$
\left(%
\begin{array}{cccc}
  0 & 0 & 1 & 0 \\
  0 & 0 & 0 & 1 \\
  1 & 0 & 1 & 0 \\
  0 & 1 & 0 & 1 \\
\end{array}%
\right)\text{ and }\left(%
\begin{array}{cccc}
  0 & 0 & 0 & 1 \\
  0 & 0 & 1 & 0 \\
  0 & 1 & 1 & 0 \\
  1 & 0 & 0 & 0 \\
\end{array}%
\right)
$$
are specialized pseudo-permutations, whereas
$$
X=\left(%
\begin{array}{ccc}
  0 & 0 & 1 \\
  0 & 1 & 0 \\
  1 & 0 & 1 \\
\end{array}%
\right)\text{ and }Y=\left(%
\begin{array}{cccc}
  0 & 0 & 0 & 1 \\
  0 & 0 & 1 & 0 \\
  0 & 1 & 1 & 0 \\
  1 & 0 & 0 & 1 \\
\end{array}%
\right)
$$
are not, in spite of being pseudo-permutations. In the first case
the index 2 is an $X$-index interior to the $X$-pair $(1,3)$; in
the second case the $Y$-pair $(2,3)$ is inside the $Y$-pair
$(1,4)$. Clearly every sub-permutation $(0,1)$-matrix is a
specialized pseudo-permutation.

\begin{thm}
\label{este} Let $X\in M$ be symmetric. Suppose $\ki=2$ and
$F=F^2$. Then

(a) $X$ is $U$-congruent to a unique specialized
pseudo-permutation matrix.

(b) $X$ is $B$-congruent to a unique specialized
pseudo-permutation $(0,1)$-matrix.
\end{thm}

\noindent{\it Proof.}  It is easy to show by induction that $X$
must be $U$-congruent to a pseudo-permutation matrix $Z$. Suppose
that $Z$ is not specialized. Then there exists a problematic
$Z$-pair $(i,j)$, having either a $Z$-pair $(k,\ell)$ inside or a
$Z$-index $s$ interior to it.

In the first case, given $0\neq a\in F$, we add $a$ times row
$\ell$ to row $j$ and then $a$ times column $\ell$ to column $j$.
This congruence transformation will replace $Z_{jj}$ by
$Z_{jj}+a^2Z_{\ell,\ell}$. It will also modify the entries
$Z_{jk}$ and $Z_{j\ell}$ on row $j$, and $Z_{kj}$ and $Z_{\ell j}$
on column $j$, into non-zero entries. As $F=F^2$ we may choose $a$
so that the $(j,j)$ entry of $Z$ becomes 0. We can then use
$Z_{ji}$ and $Z_{ij}$ to eliminate the above four spoiled entries.

In the second case we reason analogously, using the entry
$Z_{s,s}$ to eliminate the entry $Z_{jj}$, and then $Z_{ji}$ and
$Z_{ij}$ to clear the new entries in positions $(j,s)$ and $(s,j)$
back to 0.

In either case the problematic pair $(i,j)$ ceases to be a
$Z$-pair, and all other entries of $Z$ remain the same. Repeating
this process with every problematic pair produces a specialized
pseudo-permutation matrix $U$-congruent to $X$. The proves
existence in (a). The corresponding existence result in (b)
follows at once.

We are left to demonstrate the more delicate matter of uniqueness.
Let $H$ stand for either of the groups $U$ or $B$.

Let $Y$ and $Z$ be specialized pseudo-permutation matrices which
are $H$-congruent. In the case $H=B$ we further assume that $Y,Z$
are $(0,1)$-matrices. We wish to show that $Y=Z$.

Let $\hat{Y}$ be the matrix obtained from $Y$ transforming into 0
the entry $Y_{jj}$ of any $Y$-pair $(i,j)$. Clearly $\hat{Y}$ is
$H$-equivalent to $Y$ (not to be confused with $H$-congruent to
$Y$). Moreover, $\hat{Y}$ is a sub-permutation matrix (and a
$(0,1)$-matrix if $H=B$). Let $\hat{Z}$ be constructed similarly
from $Z$. All matrices $Y,\hat{Y},Z,\hat{Z}$ are $H$-equivalent,
so the uniqueness part of Theorem \ref{zero} yields that
$\hat{Y}=\hat{Z}$.

It remains to show that if $(i,j)$ is a $Y$-pair then
$Y_{jj}=Z_{jj}$, and conversely. By symmetry of $H$-congruence,
the converse is redundant.

Suppose then that $(i,j)$ is a $Y$-pair. Aiming at a
contradiction, assume that $Y_{jj}\neq Z_{jj}$ (in the case $H=B$
we are assuming that $Z_{jj}=0$, i.e. $(i,j)$ is not a $Z$-pair).

We have $A'YA=Z$ for some $A\in H$. Using the fact that $A$ is
upper triangular, for all $1\leq u,v\leq n$ we have
$$
Z_{uv}=\underset{1\leq k\leq u}\sum\;\underset{1\leq \ell\leq
v}\sum A_{ku}Y_{k\ell} A_{\ell v}.
$$
As $\ki=2$, the entry $Z_{jj}$ simplifies to
$$
Z_{jj}=\underset{1\leq k\leq j}\sum A_{kj}^2Y_{kk}.
$$
If $H=U$ then $A_{jj}=1$ and $Z_{jj}\neq Y_{jj}$. If $H=B$ then
$Y_{jj}=1$, $A_{jj}\neq 0$ and $Z_{jj}=0$. In either case, there
must exist an index $s$ such that $1\leq s<j$ and
$A_{sj}Y_{ss}\neq 0$. Choose $s$ as small as possible subject to
these conditions.

Notice that $s\neq i$, since $(i,j)$ is a $Y$-pair, which implies
that $Y_{ii}=0$.

We claim that there exists a pair $(p,q)$ such that $1\leq p\leq
j$, $1\leq q<s$, $q<i$, $p\neq q$ and $A_{pj}Y_{pq}\neq 0$.

To prove the claim we need to analyze two cases: $i<s$ or $s<i$.

Consider first the case $i<s$. Since $Y_{ss}\neq 0$, $i<s<j$ and
$(i,j)$ does not have interior $Y$-indices, there must exist an
index $t$ such that $t<s$ and $(t,s)$ is a $Y$-pair. By
transitivity, $t<j$. Now the $Y$-pair $(t,s)$ cannot be inside
$(i,j)$, so necessarily $t<i$ (*).

Since $\hat{Y}=\hat{Z}$ the only non-zero off-diagonal entry of
$Z$ in row $j$ is $Z_{ji}$, so $Z_{jt}=0$. Thus
$$
0=Z_{jt}=\underset{1\leq k\leq j}\sum\;\underset{1\leq \ell\leq
t}\sum A_{kj}Y_{k\ell} A_{\ell t}.
$$
But $A_{sj}Y_{st}A_{tt}\neq 0$, since $A_{tt}\neq 0$ is a diagonal
entry, $(t,s)$ is a $Y$-pair, and $A_{sj}\neq 0$ by the choice of
$s$. It follows that a different summand to this must be non-zero,
that is $A_{pj}Y_{pq}A_{qt}\neq 0$ for some $1\leq p\leq j$,
$1\leq q\leq t$ and $(p,q)\neq (s,t)$.

If $p=q$ then $A_{pj}Y_{pp}\neq 0$ where $p=q\leq t<s$, against
the choice of $s$. Thus $p\neq q$.

Suppose, if possible, that $q=t$. Then $p\neq s$, since $(p,q)\neq
(s,t)$. Moreover, $Y_{pt}\neq 0$. But we also have $Y_{st}\neq 0$,
with $s\neq t$. By the nature of $Y$, this can only happen if
$p=t$. But then $p=q$, which was ruled out before. It follows that
$q<t$. Since $t<s$ and $t<i$, we infer that $q<i$ and $q<s$. This
proves the claim in this case.

Consider next the case $s<i$. Since $Y_{ss}\neq 0$, either $s$ is
a $Y$-index or there is $t<s$ such that $(t,s)$ is a $Y$-pair. In
the second alternative we argue exactly as above, starting at (*)
(the fact that $t<i$ is now obtained for free, since $t<s<i$).

Suppose thus that $s$ is a $Y$-index. The only non-zero
off-diagonal entry in row $j$ of $Z$ is again $Z_{ji}$, so
$Z_{js}=0$. Thus
$$
0=Z_{js}=\underset{1\leq k\leq j}\sum\;\underset{1\leq \ell\leq
t}\sum A_{kj}Y_{k\ell} A_{\ell t}.
$$
But $A_{sj}Y_{ss}A_{ss}\neq 0$, since $A_{ss}\neq 0$ is a diagonal
entry, and the choice of $s$ ensures $A_{sj}Y_{ss}\neq 0$. As
above, there must exist $(p,q)$ such that $A_{pj}Y_{pq}Y_{qs}\neq
0$, $1\leq p\leq j$, $1\leq q\leq s$ and $(p,q)\neq (s,s)$.

If $q=s$ then $Y_{ps}\neq 0$. But $s$ is a $Y$-index, so $p=s$,
against the fact that $(p,q)\neq (s,s)$. This shows that $q<s$.
Since $s<i$, we also have $q<i$.

If $p=q$ then $A_{pj}Y_{pp}\neq 0$ with $p=q<s$, against the
choice of $s$. Therefore $p\neq q$. This proves our claim in this
final case.

The claim being settled, we choose a pair $(p,q)$ satisfying the
stated properties with $q$ as small as possible. We next produce a
another such pair with a smaller second index, yielding the
desired contradiction.

Indeed, the only non-zero off-diagonal entry in row $j$ of $Z$ is
again $Z_{ji}$ and $q<i$, so $Z_{jq}=0$. Thus
$$
0=Z_{jq}=\underset{1\leq k\leq j}\sum\;\underset{1\leq \ell\leq
q}\sum A_{kj}Y_{k\ell} A_{\ell q}.
$$
But $A_{pj}Y_{pq}A_{qq}\neq 0$ by our choice of $q$, so there
exists $(k,\ell)$ such that $A_{kj}Y_{k\ell}A_{\ell q}\neq 0$,
$1\leq k\leq j$, $1\leq \ell\leq q$ and $(k,\ell)\neq (p,q)$.
Obviously $\ell<i$ and $\ell<s$.

If $k=\ell$ then $A_{kj}Y_{kk}\neq 0$ where $k=\ell\leq q<s$,
against the choice of $s$. Thus $k\neq \ell$.

Suppose, if possible, that $\ell=q$. Then $k\neq p$, since
$(k,\ell)\neq (p,q)$. Moreover, $Y_{kq}\neq 0$. But we also have
$Y_{pq}\neq 0$, with $p\neq q$. By the nature of $Y$, this can
only happen if $k=q$. But then $k=\ell$, which was ruled out
before. It follows that $\ell<q$. This contradicts the choice of
$q$ and completes the proof. \qed

\begin{note} Every algebraic extension $F$ of the field with 2
elements satisfies $F=F^2$. The hypothesis $F=F^2$ cannot be
dropped in Theorem \ref{este}. Indeed, suppose that $z$ is not a
square in a field $F$ of characteristic 2 (e.g. $t$ is not a
square in the field $F=K(t)$, where $K$ is a field characteristic
2 and $t$ is transcendental over $K$). Then the matrix
$$
\left(%
\begin{array}{ccc}
  0 & 0 & 1 \\
  0 & 1 & 0 \\
  1 & 0 & z \\
\end{array}%
\right)
$$
is not $B$-congruent to a specialized pseudo-permutation matrix.
\end{note}

\section{Equivalence and Congruence under Parabolic Subgroups}

We fix here a standard parabolic subgroup $P$ of $G$, generated by
$B$ and a set $J$ of transpositions of the form $(i,i+1)$, $1\leq
i<n$.

Let $W$ be the group generated by $J$. Let $O_1,...,O_r$ be the
orbits of $W$ acting on $Z$. We denote by $M_i$ the largest index
in $O_i$. Note that $W$ is isomorphic to the direct product of
symmetric groups defined on the $O_i$.

Let $Y$ be a sub-permutation $(0,1)$-matrix. We let $(f,\sigma)$
stand for a $Y$-couple, where $f:\{1,...,n\}\to \{0,1\}$.

For $1\leq i,j\leq r$ we define $Y\{i,j\}$ to be equal to the
total number of indices $k$ such that $k\in S(f)$, $k\in O_i$ and
$\sigma(k)\in O_j$.

For $C\in M$ and $1\leq i,j\leq r$ we define $C[i,j]$ to be the
rank of the $M_j\times M_i$ top left sub-matrix of $C$. We also
define $C[0,j]=0$ and $C[i,0]=0$ for $1\leq i,j\leq r$.

\begin{thm}
\label{main}
 Let $P$ be a parabolic subgroup of $G$ with Weyl group
$W$. Keep the above notation. Let $C$ and $D$ be in $M$ and let
$Y$ and $Z$ be sub-permutation $(0,1)$-matrices respectively
$B$-equivalent to them. Then the following conditions are
equivalent:

(a) $C$ and $D$ are $P$-equivalent.

(b) $C[i,j]=D[i,j]$ for all $1\leq i,j\leq r$.

(c) $Y$ and $Z$ are are $W$-equivalent.

(d) $Y\{i,j\}=Z\{i,j\}$ for all $1\leq i,j\leq r$.
\end{thm}

\noindent{\it Proof.} Taking into account the type of elementary
matrices that actually belong to $P$ we see that (a) implies (b).
The equivalence of (c) and (d) is not difficult to see. Obviously
(c) implies (a). Suppose (b) holds. Since $Y$ is a a
sub-permutation $(0,1)$-matrix, it is clear that
\begin{equation}
\label{t3} Y[i,j]-Y[i-1,j]-Y[i,j-1]+Y[i-1,j-1]=Y\{i,j\},\quad
1\leq i,j\leq r.
\end{equation}
Now $Y$ and $Z$ are $P$-equivalent to $C$ and $D$, respectively,
so the equation in (b) is valid with $C$ replaced by $Y$ and $D$
replaced by $Z$. This includes also the case when $i=0$ or $j=0$.
Then (\ref{t3}) and the corresponding formula for $Z$ yield
(d).\qed

We turn our attention to $P$-congruence. Keep the notation
preceding Theorem \ref{main} but suppose now that $Y$ is symmetric
or alternating. We may assume that $\sigma$ has order 2.

Let $\sigma'$ be the permutation obtained from $\sigma$ by
eliminating all pairs $(i,j)$ in the cycle decomposition of
$\sigma$ such that either $i,j$ are in the same $W$-orbit or
$f(i)=0$ (and hence $f(j)=0$). We call $\sigma'$ the {\em reduced
permutation} associated to $Y$.

\begin{thm}
\label{main2} Suppose that $\ki\neq 2$ and $F=F^2$. Let $C$ and
$D$ be symmetric matrices and let $Y$ and $Z$ be sub-permutation
$(0,1)$-matrices respectively $B$-congruent to them. Let $\sigma'$
and $\tau'$ be the reduced permutations associated to $Y$ and $Z$,
respectively. The following conditions are equivalent:

(a) $\sigma'$ is $W$-conjugate to $\tau'$ and $Y\{i,i\}=Z\{i,i\}$
for all $1\leq i\leq r$.

(b) $C$ and $D$ are $P$-congruent.

(c) $C$ and $D$ are $P$-equivalent.

\end{thm}

\noindent{\it Proof.} It is clear that (a) implies (b) and that
(b) implies (c). Suppose (c) holds. By Theorem \ref{main}
\begin{equation}
\label{tud} Y\{i,j\}=Z\{i,j\},\quad 1\leq i,j\leq r.
\end{equation}
Writing $\sigma'$ and $\tau'$ as a product of disjoint
transpositions, condition (\ref{tud}) ensures that the number of
transpositions $(a,b)$ where $a\in O_i$, $b\in O_j$ and $i\neq j$
is the same in both $\sigma'$ and $\tau'$. For each such pair
$(a,b)$ present in $\sigma'$ and each such pair $(c,d)$ present in
$\tau'$ we let $w(a)=c$ and $w(b)=d$. Doing this over all such
pairs and all $i\neq j$ yields an injective function $w$ from a
subset of $\{1,...,n\}$ to a subset of $\{1,...,n\}$ that
preserves all $W$-orbits. We may extend $w$ to an element, still
called $w$, of $W$. This element satisfies $w\sigma'
w^{-1}=\tau'$.\qed

A reasoning similar to the above yields

\begin{thm}
\label{main2} Let $C$ and $D$ be alternating matrices and let $Y$
and $Z$ be $(1,-1)$-matrices respectively $B$-congruent to them.
Let $\sigma'$ and $\tau'$ be the reduced permutations associated
to $Y$ and $Z$, respectively. The following conditions are
equivalent:

(a) $\sigma$ is $W$-conjugate to $\tau$ and $Y\{i,i\}=Z\{i,i\}$
for all $1\leq i\leq r$.

(b) $C$ and $D$ are $P$-congruent.

(c) $C$ and $D$ are $P$-equivalent.
\end{thm}

\begin{note} The second condition in (a) is not required for
invertible matrices in either of the above two theorems.
\end{note}

\section{Number of Orbits}

Here we count the number of certain orbits under $B$-congruence.

\begin{thm} Let $C(n)$ be the number of $B$-congruence orbits of alternating matrices. Then $C(n)$ satisfies the
recursive relation
$$
C(0)=1; \quad C(1)=1;\quad C(n)=C(n-1)+(n-1)C(n-2),\; n\geq 2.
$$
\end{thm}

\noindent{\it Proof.} Suppose $Y$ is a $(1,-1)$-matrix. If column
1 of $Y$ is 0 there are $C(n-1)$ choices for $(n-1)\times (n-1)$
matrix that remains after eliminating row and column 1 of $Y$.
Otherwise there are $n-1$ choices for the position $(i,1)$, $i>1$,
of the $-1$ on column 1 of $Y$. Every choice $(i,1)$ completely
determines rows and columns 1 and $i$ of $Y$, with $C(n-2)$
choices for the $(n-2)\times (n-2)$ matrix that remains after
eliminating them.\qed

Reasoning as above, we obtain

\begin{thm} Suppose $F=F^2$ and $\ki\neq 2$. Let $D(n)$ stand for the
number of $B$-congruence orbits of symmetric matrices. Then $D(n)$
satisfies the recursive relation
$$
D(0)=1; \quad D(1)=2;\quad D(n)=2D(n-1)+(n-1)D(n-2).
$$
\end{thm}

\smallskip

\noindent{\bf REFERENCES}

\smallskip

\small

[1] Bagno E. and Cherniavski Y. \emph{Congruence B-orbits of
symmetric matrices}, preprint.

[2] Brundan J. \emph{Dense orbits and double cosets}, Algebraic
groups and their representations (Cambridge, 1997), 259--274, NATO
Adv. Sci. Inst. Ser. C Math. Phys. Sci., 517, Kluwer Acad. Publ.,
Dordrecht, 1998.

[3] Carter R. \emph{Simple Groups of Lie Type}, John Wiley $\&$
Sons, London, 1972.

[4] Duckworth W. E. \emph{A classification of certain finite
double coset collections in the classical groups}, Bull. London
Math. Soc. {\bf 36} (2004), 758--768.

[5] Ellis R., Gohberg I. and Lay D. \emph{Factorization of block
matrices}, Linear Algebra Appl. {\bf 69} (1985), 71--93.

[6] R. Gow, M. Marjoram and A. Previtali, \emph{On the irreducible
characters of a Sylow 2-subgroup of the finite symplectic group in
characteristic 2}, J. Algebra {\bf 241} (2001), 393-409.

[7] Horn R. and Johnson C. \emph{Matrix Analysis}, Cambridge
University Press, 1985.

[8] Kaplanski, I. \emph{Linear Algebra and Geometry, a Second
Course}, Allyn and Bacon, Boston, 1969.

[9] Lawther, R. \emph{Finiteness of double coset spaces}, Proc.
London Math. Soc. (3) {\bf 79} (1999), 605--625.

[10] Springer, T. A. \emph{Some results on algebraic groups with
involutions}, Algebraic groups and related topics (Kyoto/Nagoya,
1983), 525--543, Adv. Stud. Pure Math., 6, North-Holland,
Amsterdam, 1985.

\end{document}